\documentclass[a4paper]{article}

\usepackage{amsmath,amsthm,amssymb,hyperref,cleveref,enumerate}
\usepackage[utf8]{inputenc}

\theoremstyle{plain}

\newtheorem{definition}{Definition}

\newtheorem{theorem}{Theorem}

\newcommand{\RR}{\mathbb{R}}
\newcommand{\QQ}{\mathbb{Q}}
\newcommand{\ZZ}{\mathbb{Z}}
\newcommand{\NN}{\mathbb{N}}
\newcommand{\Prd}{\mathbb{P}}

\title{Fixed Point Constructions in Tilings and Cellular Automata} 

\author{
  Ilkka Törmä \\
  Department of Mathematics and Statistics \\
  University of Turku, Turku, Finland \\
  \texttt{iatorm@utu.fi} \\
  \url{http://orcid.org/0000-0001-5541-8517}
}

\begin{document}
\maketitle

\begin{abstract}
  The fixed point construction is a method for designing tile sets and cellular automata with highly nontrivial dynamical and computational properties.
  It produces an infinite hierarchy of systems where each layer simulates the next one.
  The simulations are implemented entirely by computations of Turing machines embedded in the tilings or spacetime diagrams.
  We present an overview of the construction and list its applications in the literature.

\end{abstract}

\section{Introduction}

In this short survey we present an overview of the \emph{fixed point construction} of Wang tile sets and cellular automata, also known as \emph{programmatic self-simulation}.
It is a method of defining an infinite sequence of tile sets $(T_k)_{k \geq 0}$ with the property that each $T_k$ simulates $T_{k+1}$, in the sense that valid tilings over $T_k$ have the form of infinite regular $N_k \times N_k$ grids for some constant $N_k > 1$ and each grid cell behaves like a tile of $T_{k+1}$.
The simplest variant has $T_k = T_{k+1}$ for all $k$, so that $T_0$ simulates itself.
The defining property of the technique is that the simulations are implemented almost entirely ``in software'' by computations of embedded universal Turing machines.
This provides considerable flexibility, since modifying the next tile set $T_{k+1}$ is as easy as modifying the program of the machine that runs on $T_k$.
On the flip side, the flexibility comes with the price of nonmodularity, as nontrivial constructions are hard to reuse without presenting them in detail.
The tile sets produced by the fixed point method are also inevitably massive, and it is usually impractical to determine them exactly.
Of course, this is true for most constructions in symbolic dynamics that involve simulating nontrivial computation, but the fixed point construction stands out by requiring it even for results that do not have a computational flavor, such as the existence of an aperiodic tile set.

The fixed point construction has its roots in Kurdyumov's informal article~\cite{Ku78} on probabilistic cellular automata.
Kurdyumov's ideas were implemented rigorously by Gács in his disproof of the positive rates conjecture~\cite{Ga86,Ga01}, which uses programmatic self-simulation as the basis of an extremely intricate construction.
Durand, Romashchenko and Shen isolated this part of the construction and applied it to Wang tile sets in a series of papers~\cite{DuRoSh08,DuRoSh09,DuRoSh09b,DuRoSh10} culminating in~\cite{DuRoSh12}.
A series of further applications by them and other authors has followed.
Interestingly, many of these share a common theme: a result on tiling systems, sofic shifts or cellular automata was proved in all dimensions except one or two using a geometric construction, such as Robinson tiles~\cite{Ro71} or Mozes's realization of substitutive shifts~\cite{Mo89}, and the remaining low-dimensional cases were settled by a fixed point construction.

\section{Definitions and notation}

Let $A$ be a finite alphabet and $d \geq 1$.
A \emph{pattern} is a function $P : D \to A$ with \emph{domain} $D = D(P) \subseteq \ZZ^d$.
Patterns of domain $\ZZ^d$ are called \emph{configurations}, and they form the \emph{$d$-dimensional full shift} $A^{\ZZ^d}$.
Define the \emph{shift map} $\tau : \ZZ^d \times A^{\ZZ^d} \to A^{\ZZ^d}$ by $(\tau^{\vec n} x)_{\vec v} = x_{\vec v + \vec n}$.
A set of finite patterns $F$ defines a \emph{subshift} as the subset $X_F = \{ x \in A^{\ZZ^d} \;|\; \forall P \in F, \vec v \in \ZZ^d : (\tau^{\vec v} x)|_{D(P)} \neq P \}$ where no pattern from $F$ occurs at any position.
If $F$ is finite, $X_F$ is a \emph{shift of finite type} (SFT), and if $F$ is computably enumerable, $X_F$ is an \emph{effective subshift}.
A pattern is \emph{locally valid} for $F$ if no element of $F$ occurs in it as a sub-pattern.
The set of domain-$D$ patterns occurring in a subshift $X$ is denoted $\mathcal{L}_X(D) = \{ x|_D \;|\; x \in X \}$.
If a subshift $X$ does not properly contain another nonempty subshift, it is \emph{minimal}.
If a subshift $X$ is a union of minimal subshifts, it is \emph{quasiperiodic}.\footnote{Note that an arbitrary union of minimal subshifts is generally not a subshift.}

A particular class of 2-dimensional SFTs is given by sets of \emph{Wang tiles}, which are squares with colored edges.
Formally, a Wang tile set is a subset of $C^4$ for a finite set $C$ of colors, with the four components representing the colors of the east, north, west and south edges of the square.
We forbid exactly those $1 \times 2$ and $2 \times 1$ patterns where the adjacent edges of the two tiles have different colors.
A $d$-dimensional SFT $X$ is \emph{aperiodic} if for all $x \in X$ and $\vec p \in \ZZ^d$ there exists $\vec v \in \ZZ^d$ with $x_{\vec v + \vec p} \neq x_{\vec v}$.
A set of Wang tiles is called aperiodic if it defines a nonempty aperiodic SFT.

A function $\phi : X \to Y$ between subshifts $X \subseteq A^{\ZZ^d}, Y \subseteq B^{\ZZ^d}$ is a \emph{block map} if there is a finite \emph{neighborhood} $N \subset \ZZ^d$ and \emph{local function} $\Phi : \mathcal{L}_X(N) \to B$ such that $\phi(x)_{\vec v} = \Phi((\tau^{\vec v} x)|_N)$ for all $x \in X$ and $\vec v \in \ZZ^d$.
If $X$ is an SFT and $\phi$ is surjective, then $Y$ is a \emph{sofic shift} and $X$ is an \emph{SFT cover} of $Y$.
A \emph{cellular automaton} (CA) is a block map from a full shift to itself.
A \emph{spacetime diagram} of a CA $\phi : A^{\ZZ^d} \to A^{\ZZ^d}$ is a configuration $x \in A^{\ZZ^{d+1}}$ with $x|_{\ZZ^d \times \{i+1\}} = \phi(x|_{\ZZ^d \times \{i\}})$ for each $i \in \ZZ$.
In a spacetime diagram, the $d$-dimensional slices obtained by fixing the last coordinate form a trajectory of $\phi$.

A standard reference for one-dimensional subshifts, with a short appendix on the multidimensional setting, is~\cite{LiMa95}.
See~\cite{Ka05} for a survey on the theory of CA.
Wang tiles were first defined in~\cite{Wa61}, and aperiodic sets of Wang tiles were first constructed in~\cite{Be66}.

\section{The fixed point construction}

\subsection{Tile sets}
\label{sec:tiles}

In this section we present an outline of the fixed point construction in the context of Wang tiles 
following~\cite{DuRoSh12}.
We omit much of the detail.

\begin{definition}
  Let $T$ and $S$ be two sets of Wang tiles.
  A \emph{simulation} of $S$ by $T$ with \emph{zoom factor} $N$ is an injective function $\alpha : S \to T^{N \times N}$ such that:
  \begin{itemize}
  \item
    For any $s_1, s_2 \in S$, the horizontal concatenation $s_1 s_2$ is locally valid over $T$ if and only if $\alpha(s_1) \alpha(s_2)$ is locally valid over $S$, and similarly for their vertical concatenation.
  \item
    For every valid tiling $x \in T^{\ZZ^2}$, there is a unique vector $\vec v \in [0, N-1]^2$ such that $(\tau^{\vec v + (i N, j N)}x)|_{[0, N-1]^2} \in \alpha(S)$ for all $(i,j) \in \ZZ^2$.
  \end{itemize}
\end{definition}

Suppose first that we are given the tile set $S$ and wish to define, for each large enough zoom factor $N$, a tile set $T = T(N)$ and a simulation $\alpha : S \to T^{N \times N}$.
We will implement $T$ in a way that makes the choice $S = T$ possible with relatively minor modifications.
Assume the colors of $S$-tiles are binary strings of some length $k$, that is, elements of $\{0,1\}^k$.
First, each tile $T$ has an \emph{address} $(i,j) \in [0,N-1]^2$.
Local rules ensure that its east and north neighbors have addresses $(i+1,j)$ and $(i,j+1)$ modulo $N$, respectively.
Concretely, the set of colors of $T$ is $[0,N-1] \times C$ for some auxiliary set $C$, and each tile with address $(i,j)$ has the form $((j+1, c_e), (i+1, c_n), (j, c_w), (i, c_s))$ for some $c_e, c_n, c_w, c_s \in C$.
Then in a valid tiling $x \in T^{\ZZ^2}$, 
the tiles are partitioned into a grid of $N \times N$ blocks called \emph{macrotiles} whose southwest corners are the tiles with address $(0,0)$.
The address of each tile is its relative position within the macrotile containing it.

Each tile of a macrotile $t \in T^{N \times N}$ has a specific role that depends on its address.
A fixed rectangular subset $R = [a,b] \times [c,d] \subset [0, N-1]^2$ of addresses forms the \emph{computation zone}.
It stores a simulated computation of a Turing machine $M$, which recognizes the set $S$ in the sense that it halts on input $w \in \{0,1\}^{4k}$ if and only if $w$ encodes a tile of $S$.
A reasonable choice for the parameters is $a = c = \lfloor N/3 \rfloor$ and $b = d = \lfloor 2N/3 \rfloor$, so that the dimensions of $R$ increase linearly in $N$.
The simulation of $M$ is implemented in some standard way, such as the one presented in \cite{Ro71}.
The bottom row of $R$ initializes the computation, and each row above it simulates one computation step.
If $M$ does not halt before reaching the top row of $R$, a tiling error is produced.
As long as $N$ (and thus $R$) is large enough, the simulated machine has enough time and space to verify any input that encodes a tile of $S$.

Fix an interval $I \subset [0, N-1]$ of length $k$, the number of bits in a color of $S$.
Each tile with address in $B = (\{0, N-1\} \times I) \cup (I \times \{0, N-1\})$ stores an additional bit, called a \emph{border bit}.
Their purpose is to encode the four colors of the simulated tile of $S$.
The border bits are constrained to be equal to those of the neighboring macrotiles.
For each address in $B$, we fix a contiguous path of addresses to the bottom row of the computation cone, making sure that paths originating from distinct addresses are pairwise disjoint.
Using local rules, we force all tiles with addresses on a given path to store the same bit.
In this way the border bits are routed to the computation zone, and the resulting binary word of length $4k$ forms the input to the simulated machine $M$.
If $M$ accepts this word, then it represents a tile $s \in S$, and we set $\alpha(s) = t$.
With this construction the set of locally valid macrotiles equals $\alpha(S)$, and the adjacency rules of macrotiles are identical to those of $S$.
Thus $T$ simulates $S$ with zoom factor $N$.

We would now like to choose $S = T$, so that $T = T(N)$ simulates itself.
There are a few complications that we need to handle.
First, the colors of $T$ must be binary words of constant length.
We encode one component of the address in $\lceil \log_2 N \rceil$ bits and the remaining data in another $m = \lceil \log_2 |C| \rceil$ bits.
Second, the state set and transition rules of the Turing machine $M$, which are part of the tile set $T$, depend on $S$.
We replace $M$ by a fixed universal Turing machine $M_\mathrm{U}$ that takes as input a binary word $w \in \{0,1\}^*$ and, on a separate track of the tape, a program $p \in \{0,1\}^*$, and computes $p(w)$.
The machine $M_\mathrm{U}$ should also be efficient in the sense that for any Turing machine $M'$, there is a program $p$ such that $M_\mathrm{U}(p, w) = M'(w)$ for all $w$ and the number of computation steps in $M_\mathrm{U}(p, w)$ is polynomial in $|p|$, $|w|$, and the number of computation steps in $M'(w)$.
Now the set $C$ of auxiliary colors of $T(N)$ (and thus the number $m$) is fixed, and to simulate $T(N)$ itself, it suffices to find a suitable program $p_T = p_{T(N)}$.
We enforce the simulation by requiring that each tile with address $(a+i, b)$ stores the bit $(p_T)_i$ on the program track of the simulated tape of $M_\mathrm{U}$, where $(a,b)$ is the southeast corner of the computation zone.
We call this \emph{the program condition}, and it is part of the definition of $T$.

Since we introduced the program condition into $T$, we must enforce it in the simulated version as well.
Namely, suppose for a moment that the machine $M_\mathrm{U}$ always has enough time and space to finish its simulated computation in a macrotile.
At this point of the construction we then have a tile set $T$ that simulates `tile set $T$ without the program condition,' which is of course not equal to $T$.
We can enforce the program condition in the simulated tiling as well by modifying the program $p_T$ as follows.
If the machine $M_\mathrm{U}$ simulated inside a macrotile $t \in T^{N \times N}$ is given inputs that correspond to a tile with address $(a+i, b)$ on the bottom row of the computation zone, then it reads the program bit $d \in \{0,1\}$ of any tile of $t$ with address $(a+i, b')$ for $b' \in [0, N-1]$ (which are all equal), and verifies that the simulated program bit of $t$ equals $d$.
The program condition of $T$ ensures that the initial tape of the machine $M_\mathrm{U}$ contains the program $p_T$, so we are guaranteed that $d = (p_T)_i$.
After this modification, $T$ simulates itself (given the assumption on $M_\mathrm{U}$).

Finally, consider the time and space requirements of $M_\mathrm{U}$.
The number of bits needed to store a color of $T(N)$ is $\ell = \lceil \log_2 N \rceil + m = O(\log N)$, since $m$ is fixed.
The program $p_T$ can store the number $N$ in a variable, so its length is likewise $O(\log N)$.
It should be clear that for each $w \in \{0,1\}^{4 \ell}$ the computation of $p_T(w)$ runs in time and space $O(\text{poly}(\log N))$ as long as we implement $T$ and $p_T$ in a reasonable way.
For large enough $N$ the computation zone can accommodate the computation, and then $T = T(N)$ is a self-simulating tile set.
Each valid tiling over $T$ is divided into macrotiles of level $1$, which form macrotiles of level $2$, which form macrotiles of level $3$, and so on.
In particular, the tile set is aperiodic.
In the standard terminology, a level-$k$ macrotile is a \emph{child} of the level-$(k+1)$ macrotile that contains it.

The simulation function $\alpha : T \to T^{N \times N}$ can be seen as a two-dimensional \emph{substitution}.
We can iterate it to obtain substitutions $\alpha^k : T \to T^{N^k \times N^k}$ for $k \geq 0$, and define the \emph{substitutive subshift} $X_\alpha \subset T^{\ZZ^2}$ by forbidding all finite patterns that do not occur in any $\alpha^k(t)$ for $k \geq 0$ and $t \in T$.
The set of valid tilings over $T$ contains $X_\alpha$, and $T$ can quite easily be modified to guarantee that these subshifts are equal.
In~\cite{DuRo21}, Durand and Romashchenko presented a variant of the construction in which $\alpha$ is \emph{primitive}, meaning that for some $k \geq 0$, every $t \in T$ occurs in $\alpha^k(s)$ for every $s \in T$.
In this case $X_\alpha$ is a minimal subshift.
We will discuss their results more thoroughly in \cref{sec:sofics}.

\subsection{Cellular automata}
\label{sec:ca}

The main idea of the fixed point construction for one-dimensional cellular automata is the same as for tile sets, and in fact the ``good'' spacetime diagrams of a fixed point CA have a similar substitutive structure to those of fixed point tilings, except that the width and height of the macrotiles are distinct.
The usual terminology is somewhat different, as is the implementation of information transfer between macrotiles.
As before, we first construct a CA $\phi : A^\ZZ \to A^\ZZ$ that simulates a given CA $\psi : B^\ZZ \to B^\ZZ$.
The notion of simulation has long been used informally in CA literature, but for the sake of being exact, let us use a version of \emph{injective bulking} as studied in~\cite{DeMaOlTh11}.

\begin{definition}
  A cellular automaton $\phi : A^\ZZ \to A^\ZZ$ \emph{simulates} another CA $\psi : B^\ZZ \to B^\ZZ$ if there exist integers $Q, U \geq 1$ and $s \in [0, Q-1]$ and an injection $\alpha : B \to A^Q$ such that $\phi^U(\tau^s(\alpha(x))) = \alpha(\psi(x))$ for all $x \in B^\ZZ$.
  Here $\alpha(x) \in A^\ZZ$ is defined by concatenating the blocks $\alpha(x_i)$ for $i \in \ZZ$.
\end{definition}

We assume that $\psi$ has neighborhood $N = \{-1,0,1\}$ and its state set consists of binary words of constant length $k$, and construct $\phi$ so that it also has neighborhood $N$.
The simulation parameters $Q < U$ are called the \emph{colony size} and \emph{work period}, and we choose $s = 0$.
The work period $U$ is typically vastly larger than $Q$, but it is possible to implement the CA so that $Q/U$ is arbitrarily close to 1.
Every state of $\phi$ consists of an element of $[0, Q-1] \times [0, U-1]$ called its \emph{address} and \emph{age}\footnote{Notice the difference in terminology to the tile set construction, where \emph{address} referred to a two-dimensional vector.} plus some additional data.
Say that a configuration $x \in A^\ZZ$ is \emph{valid at $i \in \ZZ$} if the states $x_{i-1} x_i x_{i+1}$ have addresses $j, j+1, j+2 \bmod Q$ for some $j$ and their ages are equal.
We only define the prototypical CA $\phi$ on configurations that are valid everywhere.
In this case, the address of a cell never changes and its age increases by $1$ modulo $U$ on every time step.

A \emph{colony} is a sequence of $Q$ adjacent cells with addresses $0, 1, 2, \ldots, Q-1$ and equal ages.
Each cell of a colony stores a \emph{simulation bit}, and the $k$ leftmost bits form the simulated state of the colony.
The ``life cycle'' of a colony, usually called a \emph{work period} when it cannot be confused with $U$, begins at age $0$ and ends at age $U-1$.
At the first step of the work period when each cell has age $0$, the simulation bits of the cells are copied onto two separate tracks called the \emph{left and right mailboxes}.
If the age of a cell is between $1$ and $Q$, it copies the contents of these tracks from its left and right neighbors.
Otherwise they remain unchanged.
In this way, during the first $Q+1$ steps of its life cycle every colony transmits its simulation bits to its nearest neighbors.
The details of this process vary between implementations.
Here we have followed~\cite{To15} and in particular assumed that $\phi$ is synchronous and deterministic.
If it is asynchronous or subject to local errors, a different scheme is needed to minimize loss of transmitted data.

Each cell of a colony also stores a \emph{program bit}, and these bits together form a program $p \in \{0,1\}^*$.
At age $Q+1$, the leftmost cell of each colony initializes the computation of the efficient universal Turing machine $M_\mathrm{U}$.
The simulated head of $M_\mathrm{U}$ is usually called the \emph{agent}.
The agent scans the program bits, simulation bits and mailboxes of the colony, computes a new simulated state $w \in \{0,1\}^k$ according to the program $p$, and writes $w$ onto the simulation bits of the $k$ leftmost cells.
We assume that the agent never steps outside the colony and halts before the work period ends at age $U-1$.
Then the age of the colony resets to $0$ and its life cycle begins anew.

As in the previous section, the CA $\phi$ simulates $\psi$ as long as $Q$ and $U$ are large enough to satisfy the time and space requirements of $M_\mathrm{U}$.
A suitable choice of the parameters allows self-simulation, meaning that we can choose $\phi = \psi$.
Note that a correct simulation hierarchy requires a program $p_\phi$ implementing $\phi$ to be stored in every colony of the initial configuration, on every simulation level.
On the first level this can be enforced by restricting the state set of $\phi$ analogously to the program condition of the tile set construction, and then $M_\mathrm{U}$ can check whether the next simulation level uses an identical program.
It depends entirely on the application how the CA behaves if the programs do not match, or if the simulation structure is invalid in some other way: it can try to correct the error, produce a special error state, or have only a partially defined local rule.

\subsection{Variable zoom factor and communication between simulation levels}

The construction is flexible enough that instead of having the tile set or the CA simulate itself, we can modify the simulated system depending on the simulation level, resulting in a sequence of tile sets $(T_k)_{k \geq 0}$ or automata $(\phi_k)_{k \geq 0}$.
The motivation is that the tile set or CA has some desirable property $P(N)$ that is limited by the zoom factor $N$, and by increasing the zoom factor on consecutive simulation levels it may be possible to obtain $P(N)$ for all $N \in \NN$.
A common use case is that $P(N)$ has a computational flavor, such as tiling the plane only if a fixed Turing machine halts within $N$ steps.

We discuss the implementation only in the case of tile sets, as the same ideas apply to the CA construction.
Our goal is then to define a sequence $(T_k, N_k)_{k \geq 0}$ of tile sets and numbers such that each tile set $T_k$ simulates $T_{k+1}$ with zoom factor $N_k$, and we would like to have as much freedom as possible in choosing the latter.
For this, we modify the construction of \cref{sec:tiles} as follows.
The program $p_T$ has an additional integer parameter: $p_T(k, w)$ determines if $w$ encodes a tile of $T_k$.
The program condition is extended so that in addition to the program $p_T$, each tile of $T_0$ stores the number $0$.
For each $k \geq 0$, the tile set $T_k$ implements a variant of the program condition where each tile of $T_{k+1}$ stores $p_T$ and the number $k+1$.
This ensures that all macrotiles have access to their simulation level.

The constraints on the zoom factors are slightly different than in the uniform case.
First, both $N_{k+1}$ and $T_{k+1}$ should be computable from $k$ in time and space $O(N_k)$.
Second, $N_k = \Omega(\log N_{k+1} + \log k)$ since the binary representation of any address in $[0,N_{k+1}-1]^2$ and the simulation level $k$ should fit in any level-$k$ macrotile.
To get a sense of the class of allowed sequences, $N_k = \lceil \log k \rceil$ for large enough $k$ grows too slowly, while $N_{k+1} = 2^{N_k}$ grows too fast.
The choices $N_k = \Theta(k^a)$ for $a \in \QQ_{> 0}$, $N_k = \Theta(2^k)$ and $N_k = \Theta(2^{2^{2^k}})$ are all acceptable.

In most applications, each macrotile stores some amount of additional information not listed above, and there are additional relations between the data stored by macrotiles of consecutive levels.
In some sense, a macrotile can ``communicate'' with its children.
The technical details vary, but the main idea is that each level-$k$ macrotile contains a special field in its state set which the simulated machine $M_\mathrm{U}$ of its parent is allowed to access.
The parent may thus read an arbitrary ``message'' from a subset of macrotiles on the bottom row of the computation zone.
The children can synchronize this message among themselves using local rules, so that it is visible to all of them or some convenient subset.
In the CA context, the agent of a colony can read and modify any data that was stored in the children during the last work period, so inter-level communication is slightly easier to implement.

\section{History and applications}

\subsection{Fault-tolerant cellular automata}
\label{sec:gacs}

The fixed point construction was first introduced in the context of \emph{probabilistic cellular automata} (PCA).
We give some basic definitions; see~\cite{MaMa14} for a more thorough review on the subject.
A $d$-dimensional PCA $\phi$ on $A^{\ZZ^d}$ consists of a neighborhood $N \subset \ZZ^d$ and a local rule $\Phi : A^N \to \Prd(A)$, where $\Prd(A)$ is the set of probability distributions over $A$.
The local rule is applied to every coordinate of a configuration $x \in A^{\ZZ^d}$ as in the case of deterministic CA, but the result is the probability distribution $\phi(x) \in \Prd(A^{\ZZ^d})$ where the value of each cell $\phi(x)_{\vec v}$ is drawn independently from $\Phi((\tau^{\vec v} x)|_N)$.
In this way $\phi$ lifts to a function $\phi : \Prd(A^{\ZZ^d}) \to \Prd(A^{\ZZ^d})$ on distributions, denoted $\mu \mapsto \mu \phi$.
Every PCA has at least one fixed point, called an \emph{invariant measure}.
If $\phi$ has a unique invariant measure $\mu$ and for every other measure $\nu$ the sequence $(\nu \phi^n)_{n \in \NN}$ converges weakly to $\mu$, then $\phi$ is \emph{ergodic}.
Intuitively, an ergodic PCA gradually forgets all details about its initial state. 
We say $\phi$ has \emph{positive rates} if the probability of $\Phi(P) = a$ is positive for all $P \in A^N$ and $a \in A$.
In other words, every local transformation occurs under $\phi$ with positive (but possibly very small) probability.
The \emph{positive rates conjecture} states that if a one-dimensional PCA has positive rates, then it must be ergodic.

In 1986 (and with an expanded paper in 2001), Gács disproved the conjecture by an intricate construction of a PCA that behaves like a deterministic CA except that at each time step, each cell has an extremely small but positive probability of ``making an error'' and ending up in an arbitrary state \cite{Ga86,Ga01}.
The design is based on a short paper of Kurdyumov from 1978 that likewise claims to refute the conjecture, but lacks detail \cite{Ku78}.

\begin{theorem}[Gács]
  There exists a nonergodic 1D PCA with positive rates.
\end{theorem}

The construction uses self-simulation to recognize and correct clusters of errors that occur at various scales.
Non-ergodicity is achieved by defining the state set as a product $\{0,1\} \times A$, initializing each cell in some state $x_i \in \{b\} \times A$ for the same bit $b$, and ``remembering'' the bit in the sense that the probability of seeing the opposite bit $1-b$ in any given cell at any given time is low.\footnote{In fact, the construction can remember an infinite sequence of bits.}
If the error rate is small enough, an argument from percolation theory shows that errors will almost surely occur in semi-isolated ``bursts'' that fit in finite (but unbounded) space-time rectangles.
At the lowest simulation level, each cell contains information about its own simulation values (address, age\ldots) and the bit to remember, as well as those of its nearest neighbors, and isolated individual errors can be fixed in one time step by a majority vote between three cells.
Larger bursts of errors are detected and repaired by higher-level colonies.
One of the main difficulties is that long sequences of errors can produce a segment of states that resembles a valid simulation structure.
Its incompatibility with the global hierarchy is only visible on its endpoints, where it cannot be determined locally which side belongs to the erroneous segment.
See Gray's condensed version of the construction for more (but not nearly all) technical details~\cite{Gr01}.
Gács and \c{C}apuni have applied the same ideas to fault-tolerant Turing machines~\cite{CaGa13}.

\subsection{Realizations by sofic shifts}
\label{sec:sofics}

In~\cite{Ho09}, Hochman proved that $d$-dimensional effectively closed subshifts can be realized as the \emph{projective sybdynamics} of $(d+2)$-dimensional sofic shifts, that is, sets of $d$-dimensional hyperplanes occurring in them.
The construction uses the older result of Mozes that two-dimensional substitutive subshifts are sofic~\cite{Mo89} to divide each configuration into regions of increasing size containing simulations of Turing machines.
The dimension was reduced from $d+2$ to $d+1$ independently by Aubrun and Sablik in~\cite{AuSa13} with Robinson tiles~\cite{Ro71}, and by Durand, Romashchenko and Shen in~\cite{DuRoSh12} using fixed point tilings.
It has since been used as a ``black box'' in numerous other constructions.
Incidentally, the realization result of Mozes was also reproved with a fixed point construction in~\cite{DuRoSh12}.

\begin{theorem}[Aubrun \& Sablik; Durand \& Romashchenko \& Shen]
  \label{thm:subdyn}
  Let $X \subseteq A^{\ZZ^d}$ be an effectively closed subshift.
  Then the higher-dimensional subshift $Y = \{y \in A^{\ZZ^{d+1}} \;|\; \exists x \in X : \forall \vec w \in \ZZ^{d+1} : y_{\vec w} = x_{\pi_d(\vec w)} \}$ is sofic, where $\pi_d : \ZZ^{d+1} \to \ZZ^d$ is defined by $\pi_d(i_1, \ldots, i_{d+1}) = (i_1, \ldots, i_d)$.
\end{theorem}

The article~\cite{DuRoSh12} is a culmination of a series of shorter papers by the same authors.
It presents the prototypical fixed point construction, collects their previous results into one place, and extends some of them.
In their proof of \cref{thm:subdyn} (in the case $d = 1$), there is a separate layer of $A$-tiles in which each vertical stripe is tiled with copies of a single $A$-symbol.
The simulations have zoom factors $N_k = Q^{c^k}$ for some large constant $Q$ and $c > 2$.
A macrotile of level $k$ with address $(x,y)$, where $y < \prod_{i=0}^k N_i$, is tasked with retrieving and storing the symbols of $f(k)$ consecutive columns of the $A$-layer situated $y$ steps away from its left border, where $f : \NN \to \NN$ is a very slowly growing function.
The definition of $N_k$ ensures that $N_{k+1} \ll \prod_{i=0}^k N_i$, so that in a level-$(k+1)$ macrotile $t$, every length-$f(k)$ word of the $A$-layer occurring within it (or right at the border between it and its neighbor) is stored in at least one child.
Using a series of inter-level messages, $t$ is in turn able to retrieve its designated length-$f(k+1)$ word.
It also computes a prefix of the sequence of forbidden patterns that defines $X$ and verifies that none of them occur in the word it retrieved.
As a direct application, the authors reprove the main result of~\cite{DuLeSh08}: there exists a tile set where the Kolmogorov complexity of every $N \times N$ pattern that occurs in a tiling is $\Omega(N)$, which is asymptotically optimal for an SFT.
Namely, we can define $X \subset \{0,1\}^\ZZ$ by forbidding all words $w$ whose Kolmogorov complexity is below $a |w| - b$ for suitable constants $a,b > 0$ and apply \cref{thm:subdyn}.

Some dynamical properties of the subshift $X$ in \cref{thm:subdyn} can be lifted to the SFT cover $Z$ of $Y$.
One of the main results of Durand and Romashchenko in~\cite{DuRo21} states that if $X$ is minimal or quasiperiodic, then the same property can be imposed on $Z$.
They modify the construction so that every macrotile contains a zone of ``diversification slots'' where all valid $2 \times 2$ patterns of tiles are manually forced to appear.
In the minimal case, there is the additional complication that if the positions of the macrotile grids and the layer of $A$-tiles in the SFT cover can be independent, the result is generally not minimal.
This is overcome by choosing a canonical configuration from each and implementing the minimal shift generated by their combination.
In the same article, the authors characterize the Turing degree spectra of quasiperiodic SFTs, which are exactly the upper closures of effectively closed subsets of $\{0,1\}^\NN$ (see~\cite[Theorem~4]{DuRo21} for details).
In particular, there exist nonempty quasiperiodic SFTs with only noncomputable configurations.

In~\cite{We17}, Westrick proves two sofic realization results that utilize the fixed point construction.
They concern two-dimensional binary subshifts where every connected component of 1-cells is a (possibly degenerate) square.

\begin{theorem}[Westrick]
  \label{thm:squares}
  Let $X_{\mathrm{sq}} \subseteq \{0,1\}^{\ZZ^2}$ be the 2D subshift where every connected component of 1-cells is either a finite square, a quarter-plane, a half-plane, or the entire $\ZZ^2$.
  \begin{enumerate}[(a)]
  \item
    Every effectively closed subshift $X \subseteq X_{\mathrm{sq}}$ that forbids all pairs of distinct squares of identical size (and possibly other patterns) is sofic.
  \item
    For $S \subseteq \NN$, let $X_S \subseteq X_{\mathrm{sq}}$ be the subshift that forbids the $n \times n$ squares for all $n \in \NN \setminus S$.
    If $S$ is computably enumerable, then $X_S$ is sofic.
  \end{enumerate}
\end{theorem}

We only sketch the construction proving the first result, as the other is quite similar.
The subshift $X_{\mathrm{sq}}$ is easily seen to be sofic, so we focus on the additional constraints of $X$.
Westrick uses the variable zoom factors $N_k = 2^{2^{2^k}}$.
In addition to the simulation structure, each macrotile $t$ stores the following information:
first, the size and relative position of each finite square of 1-cells that has at least one side completely inside $t$ or one of its neighbors;
second, the orientation and relative position of each corner and horizontal or vertical boundary of 1-cells in $t$ or its neighbors that if not part of an already stored finite square;
third, auxiliary bits that help enforce the consistency of this information between $t$ and its parent.

If the sizes of all squares are distinct, then each level-$k$ macrotile can contain at most $O(L_k^{2/3})$ squares, where $L_k = \prod_{i=0}^k N_k$ is its side length, and their sizes and positions can be encoded in $O(L_k^{2/3} \log L_k)$ bits.
This data is small enough to fit in the macrotile, but transporting it from its children to the computation zone is nontrivial, since it is too large to be simply broadcast to each child.
Westrick uses a graph theoretical argument to show that one can route the information in a nondeterministic way.
Additionally, the universal Turing machine $M_\mathrm{U}$ needs to be replaced by a multi-tape machine to speed up the computation of the relative offsets of the squares.
As in \cref{thm:subdyn}, the number of additional time steps used for computing the forbidden patterns of $X$ is increased at each simulation level.
The patterns formed by small squares are checked first, so any forbidden pattern will be eventually discovered and leads to a tiling error.

\subsection{Robust tilings}

In addition to \cref{thm:subdyn} and related realization results, the paper~\cite{DuRoSh12} presents tile sets that are \emph{robust with respect to random errors}, meaning that a configuration that is locally valid on a large subset of $\ZZ^2$ is guaranteed to be equal to a single valid tiling in a large\footnote{Large in a slightly different sense: the first subset comes from a measure-1 set with respect to a biased Bernoulli distribution, the second has asymptotic density close to 1.} subset of $\ZZ^2$.
In fact, the main result is the construction of a tile set that combines robustness with the aforementioned property that every tiling is algorithmically complex.

\begin{definition}
  Let $T$ be a set of Wang tiles and $E \subseteq \ZZ^2$.
  A $(T,E)$-tiling is a configuration $x \in T^{\ZZ^2}$ such that for all neighboring coordinates $\vec v, \vec w \in \ZZ^2 \setminus E$, the colors at the common border of $x_{\vec v}$ and $x_{\vec w}$ match.
  The Bernoulli distribution on subsets of $\ZZ^2$ with $\Prd(\vec v \in E) = \epsilon$ for each $\vec v$ is denoted $B_\epsilon$.
\end{definition}

\begin{theorem}[Durand \& Romashchenko \& Shen]
  \label{thm:bigone}
  There exists a tile set $T$ that tiles the plane and $a, b > 0$ with the following properties.
  In every valid $x \in T^{\ZZ^2}$, every $n \times n$ square has Kolmogorov complexity at least $a n - b$.
  For small enough $\epsilon > 0$, for $B_\epsilon$-almost every $E \subseteq \ZZ^2$, for each $(T,E)$-tiling $x$, the Kolmogorov complexity of $x_{[-n,n]^2}$ is $\Omega(n)$, and there is a valid tiling $y \in T^{\ZZ^2}$ with $\limsup_n |\{ \vec v \in [-n,n]^2 \;|\; x_{\vec v} \neq y_{\vec v} \}|/ (2n+1)^2 \leq 1/10$.
\end{theorem}

The construction of robust self-simulating tile sets turns out to be much simpler than the error-correcting CA discussed in \cref{sec:gacs}.
Similarly to that construction, if $\epsilon$ is small enough, then almost all error sets $E$ can be decomposed into isolated finite subsets that can be handled separately.
It is then enough to modify the prototypical fixed point construction so that every locally valid pattern shaped like the $5 \times 5$ annulus $[-2,2]^2 \setminus [-1,1]^2$ can be extended into a valid $5 \times 5$ rectangle in a unique way.
Namely, this property is inherited by macrotiles of all simulation levels, which can then correct arbitrarily large (but finite) incorrect patterns.
The additional condition on the Kolmogorov complexity of $(T,E)$-tilings requires a little more machinery, since the uniformity of the vertical columns that store the bits of the high-complexity shift $X$ can be broken by the error set $E$.
This is remedied by duplicating the information in several columns of macrotiles and periodically checking it against an error-correcting code.

\subsection{Unique ergodicity}

The article~\cite{To15} by Törmä is an application of the fixed point construction to the dynamics of cellular automata.
It concerns a weakening of the notion of \emph{nilpotency}.

\begin{definition}
  A CA $\phi : A^{\ZZ^d} \to A^{\ZZ^d}$ with a quiescent state $0 \in A$ is \emph{nilpotent} if there exists $n \in \NN$ such that $\phi^n(A^{\ZZ^d}) = \{0^{\ZZ^d}\}$.
  It is \emph{asymptotically nilpotent} if for each $x \in A^{\ZZ^d}$ we have $\phi^n(x) \to 0^{\ZZ^d}$ in the product topology as $n \to \infty$.
  It is \emph{uniquely ergodic} if it has a unique invariant measure.
\end{definition}

A nilpotent CA sends every initial configuration into the same uniform ``sink'' in a bounded number of steps.
Asymptotic nilpotency means that for every initial configuration, each cell will be in a nonzero state on only a finite (but not necessarily uniform) number of time steps.
It was shown in~\cite{GuRi08} that asymptotically nilpotent one-dimensional CA are in fact nilpotent, and the result was generalized in~\cite{Sa12,SaTo21}.
As for unique ergodicity, note that since $0$ is quiescent, the unique invariant measure cannot be any other than the Dirac point measure concentrated on $0^{\ZZ^d}$.
Then by classical results of ergodic theory, $\phi$ is uniquely ergodic if and only if for every initial configuration, the asymptotic proportion of time steps that any given cell spends in a nonzero state equals zero.
The main result of~\cite{To15} states that even one-dimensional uniquely ergodic CA are not necessarily nilpotent.

\begin{theorem}[Törmä]
  There exists a uniquely ergodic CA $\phi : A^\ZZ \to A^\ZZ$ that has a quiescent state and is not nilpotent.
\end{theorem}

This CA is defined by a fixed point construction with variable zoom factors.
The alphabet of the CA $\phi_k$ on each simulation level $k$ contains a designated ``blank'' state $B_k$, with $B_0 = 0$.
Each cell of a level-$k$ colony will spend a proportion $p_k$ of its work period simulating the state $B_k$, with $p_k \to 1$ as $k \to \infty$.
If a colony and both of its neighbors simulate the state $B_k$, it is temporarily replaced by a segment of $0$-cells.
Any local inconsistency in the simulation structure causes all nearby colonies to assume the blank state.
These properties ensure that every locally correct simulation structure is asymptotically infinitely sparse, and locally incorrect structures are quickly destroyed.
On the other hand, since a carefully chosen initial condition gives rise to an infinite nested simulation, the CA is not nilpotent.

\subsection{Realization of topological entropy}

The \emph{topological entropy} of a topological dynamical system measures its complexity and unpredictability.
In the case of subshifts and cellular automata, it can be characterized as a limit of counting patterns of increasing sizes that occur in the valid configurations and spacetime diagrams.

\begin{definition}
  The topological entropy of a subshift $X \subseteq A^{\ZZ^d}$ is the limit $h(X) = \lim_{n \to \infty} n^{-d} \log |\{ x_{[0,n-1]^d} \;|\; x \in X \}|$.
  The topological entropy of a CA $\phi$ on $A^{\ZZ^d}$ is $h(\phi) = \lim_{r \to \infty} \lim_{n \to \infty} n^{-1} \log |\{ (\phi^t(x)_{[0,r-1]^d})_{t=0}^{n-1} \;|\; x \in A^{\ZZ^d} \}|$.
\end{definition}

By classical results in symbolic dynamics (see Sections 4 and 11 of~\cite{LiMa95}), the entropy of a one-dimensional SFT is efficiently computable and its possible values are exactly the nonnegative integer multiples of logarithms of Perron numbers.
The entropies of two- and higher-dimensional SFTs were characterized by Hochman in~\cite{HoMe10} as the nonnegative right-computable numbers, that is, limits of computable decreasing sequences of rational numbers.
Independently in the preprint~\cite{GaSa17} using a variant of Robinson tiles and in~\cite{DuRo21} using fixed point tile sets, it was shown that they can be realized with transitive SFTs.

\begin{theorem}[Gangloff \& Sablik; Durand \& Romashchenko]
  The topological entropies of transitive 2D SFTs are exactly the nonnegative right-computable numbers.
\end{theorem}

The topological entropies of three- and higher-dimensional CA were characterized in~\cite{Ho09}.
Guillon and Zinoviadis handled the remaining dimensions, 1 and 2, in~\cite{GuZi13}.
Their proof uses the fixed point construction.

\begin{theorem}[Hochman; Guillon \& Zinoviadis]
  The topological entropies of one-dimensional CA are exactly the nonnegative right-computable numbers, and those of two- and higher-dimensional CA are exactly the limits of increasing computable sequences of nonnegative right-computable numbers (including $\infty$).
\end{theorem}

The high-level idea in each construction is the same.
In the SFT case, given a number $\alpha \geq 0$ in the appropriate class, construct an SFT that factors onto a subshift $X \subseteq \{0,1\}^{\ZZ^d}$ where the maximal asymptotic density of $1$-symbols in a configuration equals $\alpha$, and then allow any subset of $1$s to be independently replaced by $2$s.
In the case of cellular automata, a $d$-dimensional self-simulating CA with variable zoom factor is used to analyze a separate binary layer of the configuration, which  is not modified by the CA except when erased by the error state.
The CA produces a spreading error state if the local density of $1$-symbols is too high, which guarantees that such patterns are transient and thus do not contribute to the entropy.
This density is converted into entropy by replacing an arbitrary subset of $1$s by $2$s and composing the CA with a shift.

\subsection{Expansive directions}

Subshifts are \emph{expansive} dynamical systems: a configuration $x \in X \subseteq A^{\ZZ^d}$ can be reconstructed from any family of good enough approximations of the shifts $\tau^{\vec v} x$ for $\vec v \in \ZZ^d$.
Boyle and Lind studied a stronger variant called \emph{directional expansivity}, where a certain subset of the shifts suffices to determine $x$, in~\cite{BoLi97}.
For convenience, we define it only in the two-dimensional case.

\begin{definition}
  Let $\vec 0 \in L \subset \RR^2$ be a line through the origin with direction (slope) $\ell \in \dot \RR = \RR \cup \{\infty\}$.
  We say $\ell$ is \emph{expansive} for a subshift $X \subset A^{\ZZ^2}$, if there exists $r \geq 0$ such that the function $x \mapsto x|_D$ is injective on $X$, where $D = (L + [-r,r]^2) \cap \ZZ^2$.
  Denote the nonexpansive directions by $N(X) \subseteq \dot \RR$.
\end{definition}

For any infinite two-dimensional subshift $X$, the set $N(X)$ is nonempty and closed.
Boyle, Lind and Hochman proved in~\cite{BoLi97,Ho11} that every nonempty closed subset $L \subseteq \dot \RR$ occurs as $N(X)$ for some subshift $X$.
The subshifts constructed by Boyle and Lind consist of two layers, each of which contains regularly spaced discrete approximation of lines, and it realizes every $L$ except singleton irrational directions.
Hochman handled the remaining case with a subshift resembling the set of spacetime diagrams of a fixed point CA with variable zoom factor, but without the simulated Turing machines.
Since $X$ is not required to be of finite type, the hierarchical structure can be enforced.
Building on Hochman's proof, Guillon and Zinoviadis have characterized the sets of nonexpansive directions of SFTs in an unpublished manuscript~\cite{GuZiUP}, which also contains several realization results of SFTs with a unique nonexpansive direction.
The results are reported in Zinoviadis's PhD dissertation~\cite{Zi16}.
In the next result, a subset $A \subseteq \dot \RR$ is \emph{effectively closed} if there exists a computable sequence $(I_n)_{n \in \NN}$ of open intervals (which may have the form $(a, \infty] \cup (-\infty, b)$) with rational endpoints such that $\dot \RR \setminus A = \bigcup_{n \in \NN} I_n$.

\begin{theorem}[Guillon \& Zinoviadis]
  The sets $N(X)$ for two-dimensional SFTs $X$ are exactly the effectively closed subsets of $\dot \RR$.
\end{theorem}

The main technical contribution of the dissertation is a fixed point construction for reversible partial cellular automata (RPCA), which are defined by a local rule that is a partial function.
Note that the fixed point CA presented in \cref{sec:ca} is not reversible.
The set of spacetime diagrams of a one-dimensional RPCA is a two-dimensional SFT for which all directions $-1 < \ell < 1$ are expansive.
For an RPCA obtained by the fixed point construction with variable zoom factor, if $\prod_{k=0}^\infty Q_k/U_k = 0$, then every direction except $\infty$ (the vertical line) is expansive.
A shift map by $D_k \in \ZZ$ cells is then applied to each level $k$ of the simulation hierarchy, which changes the direction of the unique nonexpansive line.
These numbers and the simulation parameters $Q_k, U_k$ are allowed to vary within a suitable set of sequences computed by the simulated Turing machines, which provides enough control on the nonexpansive directions to realize any effectively closed set.

\bibliography{bib-fp}

\end{document}